\documentclass[11pt]{amsart}
\usepackage{amsmath,amsthm}
\usepackage{amssymb}
\usepackage[latin1]{inputenc}
\usepackage{enumerate}
\usepackage[left=3cm,top=4cm,bottom=2.5 cm,right=3cm]{geometry}
\theoremstyle{plain}
\newtheorem{theo}{Theorem}
\title{The Order of Curvature Operators on Loop Groups}
\author{Andr\'es Larra\'in-Hubach}
\address{Department of Mathematics and Statistics at Boston University
111 Cummington Street, Boston, MA 02215}
\email{alh@bu.edu}
\thanks{I would like to thank the referee for helpful suggestions.}
\begin{document}

\begin{abstract}
For loop groups (free and based), we compute the exact order of the curvature operator of the Levi-Civita connection depending on a Sobolev space parameter. This extends results of Freed \cite{F} and Maeda-Rosenberg-Torres \cite{RT}.
\end{abstract}
\maketitle

\section{Introduction}

In \cite{F}, Freed introduced a family of metrics on the free and based loop spaces of a Lie group, depending on a Sobolev space parameter $s$. These metrics determine Levi-Civita connections and curvature operators, which take values in pseudodifferential operators. In \cite{RT}, it is shown that the curvature operator is smoothing (its order is $-\infty$) when $s=1$. By \cite{F},  the order of the curvature operator is at most $\max\{-1,-2s\}$ for all $s>\frac{1}{2}$. This strange fluctuation raised the question in \cite{RT} of computing the exact order for each parameter.  Our main result (Theorem 1) is that the order of the curvature operator is at most $-2$ for $s>1$ and that this bound is exact in a specific example. Among other things, this implies that the curvature operator is trace class for $s>1$. Previously it was only known that the trace of the curvature is trace class. We also give an example where the order of the curvature operator evaluated on specific vector fields is exactly $-2s$ for $\frac{1}{2}<s<1$. Putting all these results together we get a complete answer to the original question.

Loop groups appear in several places in the mathematical physics literature. Freed's paper deals mainly with the particular values $s=\frac{1}{2}$ and $s=1$ of the Sobolev parameter. The value $s=\frac{1}{2}$, which is the critical Sobolev exponent, is important since for this metric the loop group turns out to be K\"{a}hler. The curvature associated to the Levi-Civita connection in the case $s=\frac{1}{2}$ coincides with the curvature of the canonical line bundle associated to a central extension $\widehat{LG}$ of the loop group.  It would be very interesting to understand the geometric significance of the case $s=1$, which as mentioned is special from the analytic point of view.

In \cite{PS}, the simplest Wess-Zumino-Witten model $C^\infty(S^1, G)$ is quantized using the space of sections of the determinant line bundle over an infinite Grassmanian. The  symmetry group of this line bundle corresponds to  $\widehat{LG}$. It seems important to understand the geometry of this extension for some canonical metric associated to the Sobolev metrics of $LG$. The strategy would be to relate the curvature of $\widehat{LG}$ to the curvature of the base space $LG$ using the O'Neill formulas. Therefore, this paper should be viewed as a first step towards understanding how the geometry of the symmetry group of the easiest WZW model changes with different metrics.
The necessity of quantizing higher dimensional spaces of maps $Map(\Sigma, G)$, the string theory WZW model, justifies using  the so called fractional loop groups \cite{M}. These are isometric to   the Sobolev loop groups treated here, but the interesting cases use values of the parameter $0<s<1$.

Another important application of loop groups relates to moduli spaces of flat connections of principal bundles over Riemann surfaces. Since this moduli space is the critical set of the classical Yang-Mills action, it is a fundamental object in quantum gauge theories for surfaces. In \cite{MeW}, the authors take a trivial principal $G$-bundle $\Sigma\times G\longrightarrow\Sigma$ where $\Sigma$ is a Riemann surface of genus $g$ with $b$ boundary components and $G$ is a a connected and simply-connected compact Lie group. (Note that quantization of the space of sections of this bundle is exactly the WZW model considered in \cite{M}.)
The space of connections on this bundle
is identified with $\Omega_s^1(\Sigma,\mathfrak{g})$, the space of Lie algebra valued one-forms in an $s$-Sobolev class for a sufficiently large parameter $s$. The set of flat connections $A_F(\Sigma)$ is invariant under the action of the gauge group $\mathcal{G}(\Sigma)=Map^s(\Sigma,G)$ (the gauge transformations of the same Sobolev type). For $\mathcal{G}_\partial(\Sigma)$  the kernel of the restriction to the boundary $\mathcal{G}(\Sigma)\longrightarrow\mathcal{G}(\partial\Sigma)$, the moduli space $\mathcal{M}(\Sigma)$ is defined as the quotient $A_F(\Sigma)/\mathcal{G}_\partial(\Sigma)$. When the Riemman surface has only one boundary component, this moduli space carries an action of the group $\mathcal{G}(\partial\Sigma)=L^sG$. In  analogy  to the quantization for loop groups, a central extension of $\mathcal{G}(\partial\Sigma)$  is the symmetry group of a line bundle over the moduli space $\mathcal{M}(\Sigma)$. Moreover, in the basic case where the Riemann surface is a disk, the moduli space is diffeomorphic to the based loop group $\Omega^sG$ treated here. Once again, our result on the curvature of loop groups can be seen as a first approximation to understanding the geometry of the symmetry group of these moduli spaces.

\section{Preliminaries on Loop Groups}
Following \cite{F}, we consider a completion of the space $C^\infty(S^1,G)$ of smooth loops on a connected Lie group $G$ with respect to a certain inner product given below.

 Given a loop $\gamma\in C^\infty(S^1,G)$, its formal tangent space $T_\gamma LG$ is the space of sections of the pullback bundle $\gamma^*(TG)$. Given two elements $X,Y$ of $T_\gamma LG$, an inner product is given by
\begin{align*}
\langle X,Y\rangle_s=\frac{1}{2\pi}\int_0^{2\pi}\langle(1+\Delta)^sX_t,Y_t\rangle_{\gamma(t)}dt.
\end{align*}
In the integral, we are using a left invariant inner product on $G$ which we assume is $Ad$-invariant. The operator $(1+\Delta)^s$ is defined using Seeley's theory of complex powers of elliptic operators \cite{S}. In particular, $(1+\Delta)^s$ is a pseudodifferential operator ($\psi$DO) of order $2s$. The tangent space to the loop $\gamma$, denoted $H^s_\gamma(LG)$, is the $s$-Sobolev space of sections obtained by completing  $\Gamma(T_\gamma LG)$ with respect to the inner product defined above. Here $s$ is a positive real number with $s>\frac{1}{2}$ so that, by Rellich's lemma, we get continuous sections.

Given a small neighborhood $U_e$ of the zero section in $H^s_e(LG)$, where $e$ is the identity loop, we get a local chart around $e$ in $LG$ by pointwise exponentiation
\begin{align*}
\exp_e:&U_e\longrightarrow LG,
X\mapsto (t\mapsto \exp_e X_t)
\end{align*}
Using left translation we get an atlas for the group $LG=L_sG$, the space of s-Sobolev free loops on $G$. We will be mainly concerned with the space of loops on $G$ based at the identity, denoted by $\Omega G=\Omega_sG$. For $\Omega G$, the inner product is given by
\begin{align*}
\langle X,Y\rangle_s=\frac{1}{2\pi}\int_0^{2\pi}\langle \Delta^sX_t,Y_t\rangle_{\gamma(t)}dt
\end{align*}
We use $\Delta^s$ because on $\Omega G$, $\Delta$ has trivial kernel and so has well defined complex powers. It is easy to verify that $\Omega G$ gets an atlas using the exponential map in the same way as the space $LG$. Since computations in $\Omega G$ are simpler, most of the computations of this paper will be done for $\Omega G$.

As in finite dimensions, we can use left translation to identify left invariant vector fields on $LG$ with the space $H^s(S^1,\mathfrak{g})$, the Sobolev completion of the space of smooth loops on the Lie algebra $\mathfrak{g}$ of $G$. For $\Omega G$, left invariant vector fields are identified with elements of $H^s(S^1,\mathfrak{g})$ vanishing at $0$ and $2\pi$. The operator $\Delta$ acting on left invariant vector fields will be given, using the previous identification, by $-\partial^2_\theta\otimes id$, where $id$ is the identity on the algebra $\mathfrak{g}$. The symbol of $\Delta$ at $(\theta,\xi)$ is $\xi^2 id$.

Using the inner products $\langle \cdot,\cdot \rangle_s$ on $LG$, we can construct the Levi-Civita connection by the formula
\begin{align}
2\nabla^s_XY=[X,Y]+ (1+\Delta)^{-s}[(1+\Delta)^{s}X,Y]\label{levicivdef}
+(1+\Delta)^{-s}[X,(1+\Delta)^{s}Y]
\end{align}
\cite{F}, where $X$, $Y$ are left invariant vector fields. Here we are using the usual six term formula of the Levi-Civita connection
\begin{align*}
2\langle\nabla_XY,Z\rangle&=X\langle Y,Z\rangle+Y\langle X,Z\rangle-Z\langle X,Y\rangle\\
&+\langle[X,Y],Z\rangle+\langle[Z,X],Y\rangle-\langle[Y,Z],X\rangle,
\end{align*}
where $\langle \cdot, \cdot \rangle$ is the Riemannian metric. Note that the first three terms vanish for left invariant vector fields and the last three give (\ref{levicivdef}).

The curvature operator is
\begin{align}
\Omega^s(X,Y)=\nabla^s_X\nabla^s_Y-\nabla^s_Y\nabla^s_X-\nabla^s_{[X,Y]}\label{curvdef}
\end{align}
As an operator on $Y$, $\nabla^s_X$ is a $\psi$DO of order zero, so $\Omega^s(X,Y)$ is at most of order zero. Freed \cite{F} showed that the order of this operator is at most $max\{-1,-2s\}$ for $s>\frac{1}{2}$. A natural question is to find the exact order of $\Omega^s(X,Y)$.

For the based loops case, the Levi-Civita connection is given similarly by
\begin{align*}
2\nabla^s_XY=[X,Y]+ \Delta^{-s}[\Delta^{s}X,Y]
+\Delta^{-s}[X,\Delta^{s}Y]
\end{align*}

\section{Symbol of the operators $(1+\Delta)^{\pm s}$}
Taking $s\gg 0$, we want first to compute the symbol of $(1+\Delta)^{-s}$. For this, we start with the symbol of
$(1+\Delta-\lambda)^{-1}$ where $\lambda$ is a complex parameter, since by definition \cite{S}
\begin{align}
(1+\Delta)^{-s}=\frac{i}{2\pi}\int_\Gamma \lambda^{-s}(1+\Delta-\lambda)^{-1}d\lambda.\label{compowdef}
\end{align}
By the meromorphic continuation of $(1+\Delta)^{-s}$ in \cite{S}, the result remains valid for all $s>0$. Here $\Gamma$ is a contour in the complex plane enclosing the spectrum of $1+\Delta$. From the identity
\begin{align*}
\sigma((1+\Delta-\lambda)^{-1}(1+\Delta-\lambda))=1
\end{align*}
and realizing that the total symbol of $(1+\Delta-\lambda)$ is given by $(\xi^{2}-\lambda)+1$, where $\lambda$ has weight $2$, we get
 \begin{align*}
 \sigma_{-2}((1+\Delta-\lambda)^{-1})=\frac{1}{(\xi^{2}-\lambda)}
\end{align*}
from the usual product formula for symbols (see (\ref{comppsido}) below).
Continuing this way, we see that $\sigma_{-3}((1+\Delta-\lambda)^{-1})=0$, and
\begin{align*}
\sigma_{-4}((1+\Delta-\lambda)^{-1})\sigma_{-2}(1+\Delta-\lambda)+ \sigma_{-2}((1+\Delta-\lambda)^{-1})\sigma_{0}(1+\Delta-\lambda)&=0\\
\sigma_{-4}((1+\Delta-\lambda)^{-1})(\xi^{2}-\lambda)+ (\xi^{2}-\lambda)^{-1}&=0.
\end{align*}
The symbols of order $-k$ for $k$ odd vanish and
\begin{align*}
\sigma_{-2k}((1+\Delta-\lambda)^{-1})=(-1)^{k+1}(\xi^2-\lambda)^{-k}.
\end{align*}
Plugging these into (\ref{compowdef}) we get
\begin{align*}
\sigma((1+\Delta)^{-s})=\frac{i}{2\pi}\int_\gamma \lambda^{-s}[(\xi^{2}-\lambda)^{-1}-(\xi^{2}-\lambda)^{-2}+(\xi^{2}-\lambda)^{-3}\dots]d\lambda.
\end{align*}
Integrating by parts gives
\begin{align}
\sigma((1+\Delta)^{-s})=\sum_k \binom{-s}{k}\xi^{-2s-2k}=(1+\xi^2)^{-s}
\end{align}

From $((1+\Delta)^{s})((1+\Delta)^{-s})=id$, we can compute the asymptotic expansion for $\sigma((1+\Delta)^{s})$, getting
\begin{align}
\sigma((1+\Delta)^{s})=\sum_k \binom{s}{k}\xi^{2s-2k}=(1+\xi^2)^{s}.
\end{align}

For the based loops case, the operators to be considered are $\Delta^{\pm s}$ and reasoning as above gives $\sigma(\Delta^{\pm s})(\theta, \xi)=(\xi^2)^{\pm \frac{s}{2}}$.

\section{The Zero-Order Terms of the Connection and Curvature}

For the Levi-Civita connection $\nabla=\nabla^s$ on $LG$, recall that
\begin{align}
2\nabla_XY=[X,Y]- (1+\Delta)^{-s}[(1+\Delta)^{s}X,Y]+(1+\Delta)^{-s}[X,(1+\Delta)^{s}Y]\label{levcivdefbasloop}
\end{align}

We label the terms on the right hand side of the equation as (a), (b), (c). The term (b) is a pseudodifferential
operator of order $-2s$, so it does not contribute to the zero order symbol. Letting ${e_i}$ be an orthonormal
basis for the Lie algebra $\mathfrak{g}$ with structure constants $\{C^k_{ij}\}$, we get
\begin{align}
\sigma_0 (\nabla_X)^i_k=C^i_{jk}X^j
\end{align}
by summing the terms corresponding to (a) and (c).

For the curvature form, we get
\begin{align}
\sigma_0(\Omega(X,Y))^i_k= \sigma_0(\nabla_X)^i_l\sigma_0(\nabla_Y)^l_k- \sigma_0(\nabla_Y)^i_l\sigma_0(\nabla_X)^l_k
-\sigma_0(\nabla_{[X,Y]})^i_k
\end{align}

This gives
\begin{align*}
\sigma_0(\Omega(X,Y))^i_k&=C^i_{jl}C^l_{mk}X^jY^m-C^i_{jl}C^l_{mk}Y^jX^m-C^i_{jk}C^j_{ml}X^mY^l\\
&=(C^i_{mj}C^j_{lk}-C^i_{lj}C^j_{mk}-C^i_{jk}C^j_{ml})X^lY^m\\
&=0
\end{align*}
by the Jacobi identity. Therefore the zero-order term of the curvature vanishes.

\section{The $(-1)$-Order Terms of the Connection and Curvature}

For the $(-1)$-order term of the connection, we only need to use the term (c). Recall the product formula
for the asymptotic expansion of the symbol of a composition of pseudodifferential operators:
\begin{align}
\sigma(P\circ Q)\sim \sum_\alpha \frac{1}{i^{|\alpha|}\alpha !}\,\partial ^\alpha _\xi \sigma(P)\cdot\partial^\alpha_\theta \sigma(Q)\label{comppsido}
\end{align}
Applying this formula to $(3)$ we get the following expansion
\begin{align*}
\sigma (\frac{1}{2}(1+\Delta)^{-s}[X,(1+\Delta)^{s}(\cdot)]) = &\frac{1}{2}(1+\xi^2)^{-s}C^i_{jk}X^j(1+\xi^2)^s\\
&-\frac{i}{2}\partial_\xi(1+\xi^2)^{-s}\partial_\theta
(C^i_{jk}X^j(1+\xi^2)^s)+\ldots
\end{align*}
For the $(-1)$-order symbol we only need the term which involves the first derivatives in both variables. We obtain
\begin{align}
\sigma_{-1}(\nabla_X)^i _k=is C^i_{jk}\dot{X}^j\xi |\xi|^{-2},
\end{align}
where $\dot{X}^j=\partial_\theta X^j$. To find the $(-1)$-order symbol of the curvature we must treat each term of (\ref{curvdef}) separately, e.g.
\begin{align*}
\sigma_{-1}(\nabla_X\nabla_Y)^i _k&=\sigma_{-1}(\nabla_X)^i _l \sigma_{0}(\nabla_Y)^l _k+\sigma_{0}(\nabla_X)^i _l\sigma_{-1}(\nabla_Y)^l _k\\
&= is C^i_{jl}\dot{X}^j\xi |\xi|^{-2}C^l_{mk}Y^m  +is C^i_{jl}X^j C^l_{mk}\dot{Y}^m\xi |\xi|^{-2}.
\end{align*}
One easily obtains
\begin{align*}
\sigma_{-1}\Omega(X,Y)^i _k&= \sigma_{-1}(\nabla_X\nabla_Y)^i _k-\sigma_{-1}(\nabla_Y\nabla_X)^i _k -\sigma_{-1}([X,Y])^i _k\\
&= is C^i_{jl}\dot{X}^j\xi |\xi|^{-2}C^l_{mk}Y^m  + is C^i_{jl}X^j C^l_{mk}\dot{Y}^m\xi |\xi|^{-2}\\
&\,\,\,\,\,\,\,- is C^i_{jl}\dot{Y}^j\xi |\xi|^{-2}C^l_{mk}X^m  -is C^i_{jl}Y^j C^l_{mk}\dot{X}^m\xi |\xi|^{-2}\\
&\,\,\,\,\,\,\,-is C^i_{jk}([X,Y]^j)^\cdot \xi |\xi|^{-2}.
\end{align*}
The last term is
\begin{align*}
-is C^i_{jk}([X,Y]^j)^\cdot \xi |\xi|^{-2}=-is C^i_{jk}C^j_{lm}(\dot{X}^lY^m+X^l\dot{Y}^m)\xi |\xi|^{-2}
\end{align*}
Combining the terms with $\dot{X}$ gives
\begin{align*}
is C^i_{jl}\dot{X}^j\xi |\xi|^{-2}C^l_{mk}Y^m -is C^i_{jl}Y^j C^l_{mk}\dot{X}^m\xi |\xi|^{-2}
-is C^i_{jk}C^j_{lm}\dot{X}^lY^m\xi |\xi|^{-2}.
\end{align*}
This can be rewritten as
\begin{align*}
is\xi |\xi|^{-2}(C^i_{jl}\dot{X}^jC^l_{mk}Y^m- C^i_{jl}Y^j C^l_{mk}\dot{X}^m-C^i_{jk}C^j_{lm}\dot{X}^lY^m),
\end{align*}
which in turn equals
\begin{align*}
is \xi |\xi|^{-2}(C^i_{jl}C^l_{mk}-C^i_{ml}C^l_{jk}-C^i_{lk}C^l_{jm})\dot{X}^jY^m=0
\end{align*}
by the Jacobi identity. Similarly, the terms with $\dot{Y}$ also vanish. This proves that the $(-1)$-order symbol of the curvature operator also vanishes.

\section{The $(-2)$-Order Symbol of the Connection and Curvature}
In this section, we consider $s> 1$, so that (b) in (\ref{levcivdefbasloop}) does not contribute to the
the $(-2)$-order symbol. For (c) in (\ref{levcivdefbasloop}), we have
\begin{align}
\sigma (\frac{1}{2}(1+\Delta)^{-s}[X,(1+\Delta)^{s}(\cdot)])=&\frac{1}{2}(1+\xi^2)^{-s}C^i_{jk}X^j(1+\xi^2)^s\notag\\
&-\frac{i}{2}\partial_\xi(1+\xi^2)^{-s}\partial_\theta
(C^i_{jk}X^j(1+\xi^2)^s)\label{asympexp}\\
&-\frac{1}{4}\partial ^2 _\xi(1+\xi^2)^{-s}\partial^2_\theta
(C^i_{jk}X^j(1+\xi^2)^s)+\ldots\notag
\end{align}
We note that
\begin{align}
(1+\xi^2)^{-s}&=\xi^{-2s}+\binom{-s}{1}\xi^{-2s-2}+\binom{-s}{2}\xi^{-2s-4}+\ldots\notag\\
\partial_\xi(1+\xi^2)^{-s}&=(-2s)\xi^{-2s-1}+\binom{-s}{1}(-2s-2)\xi^{-2s-3}+\binom{-s}{2}(-2s-4)\xi^{-2s-5}\ldots\label{asympexp2}\\
\partial ^2 _\xi(1+\xi^2)^{-s}&= (-2s)(-2s-1)\xi^{-2s-2}+\binom{-s}{1}(-2s-2)(-2s-3)\xi^{-2s-4}\ldots\notag\\
(1+\xi^2)^{s}&=\xi^{2s}+\binom{s}{1}\xi^{2s-2}+\binom{s}{2}\xi^{2s-4}+\ldots\notag
\end{align}
The terms with $\xi^{-2}$ in (a) in (\ref{levcivdefbasloop}) disappear, and using (\ref{asympexp2}) in (\ref{asympexp}) gives
\begin{align}
\sigma_{-2}(\nabla_X)^i_k=-s(s+\frac{1}{2})C^i_{jk}\ddot{X}^j\xi^{-2}.
\end{align}
We now compute $\sigma_{-2}$ for the curvature operator.
\begin{align*}
\sigma_{-2}(\nabla_X\nabla_Y)^i _k&=\sigma_{-2}(\nabla_X)^i _l \sigma_{0}(\nabla_Y)^l _k +
\sigma_{0}(\nabla_X)^i _l\sigma_{-2}(\nabla_Y)^l _k + \sigma_{-1}(\nabla_X)^i _l \sigma_{-1}(\nabla_Y)^l _k\\
&\quad +i^{-1}\partial_\xi(\sigma_{-1}(\nabla_X)^i _l)\partial_\theta( \sigma_{0}(\nabla_Y)^l _k)\\
&=-s(s+\frac{1}{2})C^i_{jl}C^l_{mk}(\ddot{X}^j\xi^{-2}Y^m+X^j\ddot{Y}^m)\xi^{-2}-s(s+1)C^i_{jl}C^l_{mk}\dot{X}^j\dot{Y}^m\xi^{-2}.\\
\end{align*}
The expression for $\sigma_{-2}(\nabla_Y\nabla_X)^i _k$ is the same as above with $X$ and $Y$ exchanged. The third term in
the curvature expression is
\begin{align*}
\sigma_{-2}(\nabla_{[X,Y]})^i_k&=-s(s+\frac{1}{2})C^i_{jk}([X,Y]^j)^{\cdot\cdot}\xi^{-2}\\
&=-s(s+\frac{1}{2})C^i_{jk}
C^j_{lm}(\ddot{X}^lY^m+2\dot{X}^l\dot{Y}^m+X^l\ddot{Y}^m)\xi^{-2}
\end{align*}
Gathering together all the terms with two derivatives of $X$, we get
\begin{align*}
-s(s+\frac{1}{2})C^i_{jl}\ddot{X}^j\xi^{-2}C^l_{mk}Y^m+s(s+\frac{1}{2})C^i_{jl}Y^jC^l_{mk}\ddot{X}^m\xi^{-2}
+s(s+\frac{1}{2})C^i_{jk}
C^j_{lm}\ddot{X}^lY^m\xi^{-2},
\end{align*}
which equals
\begin{align*}
-s(s+\frac{1}{2})\xi^{-2}(C^i_{lj}C^j_{mk}-C^i_{mj}C^j_{lk}-C^i_{jk}C^j_{lm})\ddot{X}^lY^m
\end{align*}
Once again, the expression in the middle is the Jacobi identity, so this term vanishes. Similarly, the term with two derivatives
of $Y$ is zero. However, the terms with $\dot{X}$ and $\dot{Y}$ do not vanish.
 We get
\begin{align*}
-s(s+1) C^i_{jl}\dot{X}^j C^l_{mk}\dot{Y}^m\xi^{-2}+s(s+1) C^i_{jl}\dot{Y}^j C^l_{mk}\dot{X}^m\xi^{-2}
+s(s+\frac{1}{2})C^i_{jk}
C^j_{lm}2\dot{X}^l\dot{Y}^m\xi^{-2},
\end{align*}
which can be written
\begin{align*}
-s(s+1)\xi^{-2}(C^i_{jl}\dot{X}^j C^l_{mk}\dot{Y}^m-C^i_{jl}\dot{Y}^j C^l_{mk}\dot{X}^m-C^i_{jk}
C^j_{lm}\dot{X}^l\dot{Y}^m) +s^2C^i_{jk}
C^j_{lm}\dot{X}^l\dot{Y}^m\xi^{-2}
\end{align*}
The first term again vanishes by the Jacobi identity, so finally
\begin{align}
\sigma_{-2}(\Omega(X,Y)^i_k)=s^2C^i_{jk}
C^j_{lm}\dot{X}^l\dot{Y}^m\xi^{-2}
\end{align}
We now compute this expression in a particular case, to show that the order of the curvature can be exactly $-2$. The Lie algebra $\mathfrak{su}(2)$ of the Lie group $SU(2)$ is given by three generators $e$, $f$ and $g$, subject to the relations $[e,f]=-2g$, $[e,g]=2f$ and $[f,g]=-2e$. Consider the left invariant vector fields on $LG$ given by $X=(\sin \theta)e$ and $Y=(\sin \theta)f$. With these vector fields, we get
\begin{align*}
\sigma_{-2}(\Omega(X,Y)_1^2)=(4s^2\cos^2\theta)\xi^{-2}.
\end{align*}

\section{The Based Loops Case}
For based loops, we will get analogous results.

 The zero order symbol of the connection is
given by
\begin{align*}
\sigma_0 (\nabla_X)^i_k=C^i_{jk}X^j,
\end{align*}
and the zero order term of the curvature operator again vanishes.

For the $-1$ order symbol of the connection, we only need the term (c) in (\ref{levcivdefbasloop}). We have
\begin{align*}
\frac{1}{2}\Delta^{-s}[X,\Delta^{s}Y]=\frac{1}{2}(\Delta^{-s}\circ ad_X\circ\Delta^{s})(Y).
\end{align*}
For $Q=ad_X\circ\Delta^{s}$, we have
\begin{align*}
\sigma(Q)_k^i=C^i_{jk}X^j\xi^{2s}.
\end{align*}
There are no more terms in the expansion since the symbol of $ad_X$ is independent of $\xi$.

We find that the total symbol of (c) is
\begin{align*}
\sigma(c)_k^i=\frac{1}{2}(C^i_{jk}X^j+isC^i_{jk}\dot{X}^j\xi^{-1}+s(-2s-1)C^i_{jk}\ddot{X}^j\xi^{-2}+\ldots)
\end{align*}

From this we see that the terms of order $-1$ and $-2$ of the symbol of the connection are the same as the
free loop case, so the same holds for the curvature. We conclude that, in this case also, the order of the
curvature operator is at most $-2$, and the same example as before implies that this bound is sharp. We also get the same formula as in the free loop case for the $-2$ order term. From now on, we will work only with based loops. For free loops we get the same results.

\section{The value $s=1$}
The formulas in the previous sections are valid for $s> 1$, since in this case the second term
in the right-hand side of the formula
\begin{align}
2\nabla_XY=[X,Y]- \Delta^{-s}[\Delta^{s}X,Y]+\Delta^{-s}[X,\Delta^{s}Y]\label{bslevciv}
\end{align}
is a pseudodifferential operator of order $-2s< -2$.Thus, this term does not contribute to the terms of order $-1$ and $-2$ in
the asymptotic expansion of the curvature operator. However, if we take $s=1$, we get another contribution to the term of order $-2$ in the asymptotic expansion of the
connection operator and hence, of the curvature operator.

The asymptotic expansion of the second term is given by
\begin{align}
-\sigma(\frac{1}{2}\Delta^{-1}[\Delta X,\cdot])^i_k&=-\dfrac{1}{2}(\sum_\alpha\frac{1}{i^{\mid\alpha\mid}\alpha!}\partial_\xi^\alpha\xi^{-2}
\partial_\theta^\alpha(-C^i_{jk}\ddot{X}^j))\notag\\
&=\frac{1}{2}(C^i_{jk}\ddot{X}^j\xi^{-2}+\frac{1}{i}(-2)\xi^{-3}C^i_{jk}\partial_\theta \ddot{X}^j+\ldots)\label{asymseq1}
\end{align}
This implies that the $-2$-order term of the connection equals the term we had previously (evaluated at $s=1$) plus the first term in the last line of (\ref{asymseq1}). We get in the end
\begin{align*}
\sigma_{-2}(\nabla_X)^i_k&=-\frac{3}{2}C^i_{jk}\ddot{X}^j\xi^{-2}+\frac{1}{2}C^i_{jk}\ddot{X}^j\xi^{-2}\\
&=-C^i_{jk}\ddot{X}^j\xi^{-2}
\end{align*}

Substituting this in the formula for the curvature gives
\begin{align*}
\sigma_{-2}(\nabla_X\nabla_Y)^i _k&=\sigma_{-2}(\nabla_X)^i _l \sigma_{0}(\nabla_Y)^l _k\\
&\quad +\sigma_{0}(\nabla_X)^i _l\sigma_{-2}(\nabla_Y)^l _k+ \sigma_{-1}(\nabla_X)^i _l \sigma_{-1}(\nabla_Y)^l _k\\
&\quad +i^{-1}\partial_\xi(\sigma_{-1}(\nabla_X)^i _l)\partial_\theta( \sigma_{0}(\nabla_Y)^l _k)\\
&=-C^i_{jl}\ddot{X}^j\xi^{-2}C^l_{mk}Y^m-C^i_{jl}X^jC^l_{mk}\ddot{Y}^m\xi^{-2}\\
&\quad -2C^i_{jl}\dot{X}^j C^l_{mk}\dot{Y}^m\xi^{-2}
\end{align*}
For the second term of the curvature we get the same expression with $X$ and $Y$ exchanged. Now
\begin{align*}
\sigma_{-2}(\nabla_{[X,Y]})^i_k=-C^i_{jk}([X,Y]^j)^{\cdot\cdot}\xi^{-2}=-C^i_{jk}
C^j_{lm}(\ddot{X}^lY^m+2\dot{X}^l\dot{Y}^m+X^l\ddot{Y}^m)\xi^{-2}
\end{align*}
Terms with second derivatives vanish as before while the terms with single derivatives add up to the following expression
\begin{align*}
- 2C^i_{jl}\dot{X}^j C^l_{mk}\dot{Y}^m\xi^{-2}+ 2C^i_{jl}\dot{Y}^j C^l_{mk}\dot{X}^m\xi^{-2}+2C^i_{jk}
C^j_{lm}\dot{X}^l\dot{Y}^m\xi^{-2}=0,
\end{align*}
so we can conclude that
\begin{align*}
\sigma_{-2}(\Omega(X,Y))^i_k=0.
\end{align*}
From this we see that the second order term of the curvature operator vanishes for $s=1$. This agrees with \cite{RT}, where it is proven that the curvature operator is smoothing for $s=1$.

\section{Computations for $\frac{1}{2}<s<1$}

The asymptotic expansion of the symbol of the second term in (\ref{bslevciv}) is
\begin{align*}
-\frac{1}{2}C_{jk}^i[(\Delta^{s}X^j)\xi^{-2s}-i\partial_\theta\Delta^sX^j\xi^{-2s-1}]+O(\xi^{-2s-2})
\end{align*}
The third term in (\ref{bslevciv}) is
\begin{align*}
\frac{1}{2}C_{jk}^i[X^j + i2s\dot{X}^j\xi^{-1}-s(2s+1)\ddot{X}^j\xi^{-2}]+O(\xi^{-3})
\end{align*}

The $-2s$-order term of the connection is $\sigma_{-2s}(\nabla_X)_k^i=-\frac{1}{2}C_{jk}^i(\Delta^sX^j)\xi^{-2s}$. For the curvature
 operator, the term of order zero vanishes as expected. The term of order $-2s$ is
 \begin{align*}
\sigma_{-2s}(\Omega(X,Y))^i_k&=-\frac{1}{2}C_{jl}^i(\Delta^sX^j)\xi^{-2s}C_{km}^lY^m-\frac{1}{2}C^i_{jl}X^jC_{mk}^l(\Delta^sY^m)\xi^{-2s}\\
&\quad+\frac{1}{2}C_{jl}^i(\Delta^sY^j)\xi^{-2s}C_{km}^lX^m+\frac{1}{2}C^i_{jl}Y^jC_{mk}^l(\Delta^sX^m)\xi^{-2s}\\
&\quad+\frac{1}{2}C_{jk}^i(\Delta^s[X,Y]^j)\xi^{-2s}
\end{align*}
For this expression to vanish we would need a Leibniz rule for the operator $\Delta^s$, but no such rule exists. Therefore the
term of order $-2s$ does not vanish in general and this gives the order of the curvature operator. The same vector fields used in section $6$ have nonzero curvature, so the order is exactly $-2s$ in this case.

\section{Conclusion}
We finish with the following result which extends those obtained in \cite{F} and \cite{RT}.
\begin{theo}
 The order of $\Omega^s(X,Y)$, for either based or free loop groups, is at most $-2$ for $s> 1$ and at most $-2s$ for $\frac{1}{2}< s< 1$. These bounds are optimal. For $s=1$ the order
  of the curvature is $-\infty$ for based loops.  For $s>1$, the curvature operator is trace class.
\end{theo}

\end{document}